\documentclass[12pt]{amsart} 
\usepackage{amsmath} 
\usepackage{amssymb}
\theoremstyle{plain} 
\newtheorem{theor}{Theorem}[section]
{}
\newtheorem{lemma}{Lemma}
\newtheorem{propo}{Proposition}
 
{}
\theoremstyle{definition} 
\newtheorem{defin}{Definition}[section]

\theoremstyle{remark}

\numberwithin{equation}{section}
\DeclareMathOperator{\cov}{cov} 
 
\DeclareMathOperator{\dom}{domain} \DeclareMathOperator{\cof}{cof}

\newcommand{\IF}{\text{ if }} 
\newcommand{\AND}{\text{ and }}  \newcommand{\forces}[2]{\Vdash_{#1} \mbox{``} #2 \mbox{''}}

\newcommand{\Poset}{{\mathbb P}}

\newcommand{\Qposet}{{\mathbb Q}} 
\newcommand{\Sposet}{{\mathbb S}}

\newcommand{\presup}[2]{\, ^{#1} \! #2} 

\newcommand{\fomom}{\presup{\omega}{\omega}}
\newcommand{\wfomom}{\presup{\stackrel{\omega}{\smile}}{\omega}}
  
\newcommand{\PF}{{\mathbb P}{\mathbb F}}
\thanks{The research for this paper was completed while the second
author was visiting the first author at Rutgers University funded by
NSF grant DMS97-04477. The second author's research is also funded by NSERC.
This paper is number B665 in the first author's personal listing.}
\title[Covering Numbers of Mycielski Ideals]{The Covering Numbers of
Mycielski Ideals are all Equal}
\author[S. Shelah]{Saharon Shelah} \address{
Institute of Mathematics, The Hebrew University, Jerusalem 91904, Israel
}
\email{shelah@math.huji.ac.il}
\author[J. Stepr\={a}ns]{Juris Stepr\={a}ns}
\address{Department of Mathematics, York University\\
4700 Keele Street,
Toronto, Ontario\\ Canada \ \ \  M3J 1P3}
\email{steprans@mathstat.yorku.ca}
\date{}
\begin{document}
\begin{abstract}
The Mycielski  ideal ${\mathfrak M}_k$ is defined to consist of all sets
$A\subseteq \presup{\omega}{k}$ such that
$\{f\restriction X: f \in A\} \neq
\presup{X}{k}$
 for all
$X\in[\omega]^{\aleph_0}$. 
It will be shown that the covering numbers for these ideals are all
equal. However, the covering numbers of the closely associated Ros\l
anowski ideals will be shown to be consistently different.
\end{abstract}
\maketitle
 \bibliographystyle{plain}
\section{Introduction}
In \cite{MR39:2934} J.~Mycielski defined a class of ideals which have
been studied in various contexts by several authors \cite{MR96k:03109,
MR96g:20001,
MR95m:04003,
MR95d:54033,
MR94i:03094,
MR94g:26004,
MR92e:04004,
MR91d:04001,
MR90m:28004,
MR94d:03098}. This paper is devoted to examining the covering numbers
of these ideals as well as those of a closely related class of
ideals. It will be shown that, while the covering number of the
Mycielski ideals is independent of their dimension, the covering
number of the related ideals is very closely related to their dimension.
\begin{defin}
The Mycielski \label{d:m} ideal ${\mathfrak M}_k$ is defined to
consist of all sets 
$A\subseteq \presup{\omega}{k}$ such that for all
$X\in[\omega]^{\aleph_0}$
\begin{equation}\label{e:1}\{f\restriction X: f \in A\} \neq
\presup{X}{k}
\end{equation} A function $\Phi$ on $[\omega]^{\aleph_0}$ will be said
to witness that $A \in {\mathfrak M}_k$ if $\Phi(X) \in
\presup{X}{k} \setminus \{f\restriction X: f \in A\}$ for each $X \in
[\omega]^{\aleph_0}$.
\end{defin}
Notice that if $A \in {\mathfrak M}_k$ and $X$ is an infinite subset
of $\omega$ then not only is there some $g \in \presup{X}{k}$ such
that for all $f \in A$ there is some $x\in X$ such that $f(x) \neq
g(x)$ but, in fact, there some $g \in \presup{X}{k}$ such
that for all $f \in A$ there are infinitely many $x\in X$ such that $f(x) \neq
g(x)$. The next definition will generalize this version of the
Mycielski ideals. 

\begin{defin}
 Let $\PF_k$ denote the set of all \label{d:r} functions $f:X \to k$
where $X$ is 
a coinfinite subset of $\omega$.
The Ros\l onowski ideal ${\mathfrak R}_k$ is defined to consist of all sets
$A\subseteq \presup{\omega}{k}$ such that for all
 $g\in\PF_k$ there is
an extension $g'\supseteq^* g$ such that $g'\in\PF_k$ and 
$g'\not\subseteq^* f$ for all $f \in A$.
 A function $\Phi$ on $\PF_k$ will be said
to witness that $A \in R(k)$ if $g\subseteq \Phi(g) \in
\PF_{k} $ for each $g \in \PF_k$ and $\Phi(g) \not\subseteq^* f$ for all
$f \in A$.
\end{defin}
It is worth noting that neither of these ideals has a simple
definition. Indeed, since the definition given is $\Pi^1_2$ many of
the usual arguments which apply to Borel ideals must be applied with
great care, if at all, in this context. For an alternate approach to
finding a nice base for the Mycielski ideals see \cite{MR95d:54033}.

The covering numbers of the ideals ${\mathfrak R}_k$ have a connection
to  gaps in ${\mathcal P}(\omega)/[\omega]^{<\aleph_0}$. Indeed, the
assertion that $\cov({\mathfrak R}_2)= \aleph_1$ can be interpreted as saying
there are many Hausdorff gaps. To see this,
suppose that $\{A_\xi\}_{\xi\in\omega_1}$ is a cover of $2^\omega$ by
sets in ${\mathfrak R}_2$ witnessed by
$\{\Phi_\xi\}_{\xi\in\omega_1}$. If $\{f_\xi\}_{\xi\in\omega_1}$
 is any $\subseteq^*$-increasing sequence  in $\PF_k$ such that
$f_{\xi+1} = \Phi_\xi(f_xi)$ then 
$\{(f_\xi^{-1}\{0\},f_\xi^{-1}\{1\})\}_{\xi\in\omega_1}$ is a
Hausdorff gap. Hence a large tree all of whose branches are Hausdorff
gaps can be constructed using $\cov({\mathfrak R}_2)=
\aleph_1$. It will be shown that similar assertions for
$\cov({\mathfrak R}_n)= \aleph_1$ are not equivalent to
$\cov({\mathfrak R}_2)= \aleph_1$ for $n> 2$. 
\section{Equality and inequality}
\begin{theor}
If $m$ and $n$ are integers greater than 1 then
 $\cov({\mathfrak M}_k) = \cov({\mathfrak M}_n)$.\end{theor}
\begin{proof} To begin, notice that if $\Phi$ witnesses that $A\in
{\mathfrak M}_k$ then 
$$\{f\in \presup{\omega}{k+1}: (\forall X \in [\omega]^{\aleph_0})
f\restriction X \neq \Phi(X)\}$$ belongs to ${\mathfrak M}_{k+1}$. It follows
that $\cov({\mathfrak M}_k) \geq \cov({\mathfrak M}_{k+1})$. It therefore suffices to show that
$\cov({\mathfrak M}_{k^2}) \geq \cov({\mathfrak M}_k)$ for each $n\geq 2$.

To this end, let $\beta : \omega \to [\omega]^2$ be a bijection and
let $\beta_s(n)$ be the smallest member of $\beta(n)$ and
$\beta_g(n)$ be the greatest member of $\beta(n)$. Define a relation
$\equiv_\beta$ on partial functions from $\omega$ to $k$ and partial functions
from $\omega$ to $k^2$ by $f\equiv_\beta g$ if and only if the following
conditions \eqref{c:1} and \eqref{c:2} hold:
\begin{equation}\label{c:1}
  (\forall \{n,m\}\in [\dom(g)]^2) \beta(n) \cap \beta(m) = \emptyset
\end{equation}
\begin{equation}\label{c:2}
  (\forall n \in \dom(g)) g(n) = kf(\beta_s(n)) + f(\beta_g(n))
\end{equation}
Now suppose that $\mathcal A$ is a cover of $\presup{\omega}{(k^2)}$ by
sets in ${\mathfrak M}_{k^2}$ and that $\Phi_A$ witnesses that $A\in {\mathfrak M}_{k^2}$ for
each $A\in \mathcal A$. Now, for $A\in \mathcal A$ define 
\begin{equation}\label{e:star}
A^* = \{f \in \presup{\omega}{k} : (\forall X \in [\omega]^{\aleph_0})
(\forall Z \in [\omega]^{\aleph_0}) f\restriction X\not\equiv_\beta
\Phi(Z)\}
\end{equation}
It will be shown that $\{A^* : A\in {\mathcal A}\}$ is a cover of
$\presup{\omega}{k}$ by sets in the ideal ${\mathfrak M}_k$. 

To see that each $A^* \in {\mathfrak M}_k$ let $A\in {\mathfrak M}_{k^2}$ and $X\in
[\omega]^{\aleph_0}$. Let $\{\{e_i, d_i\}\}_{i\in\omega}$ be disjoint
pairs from $X$ such that $e_i < d_i$ for all $i$. Let $Z =
\{\beta^{-1}(\{e_i,d_i\})\}_{i\in\omega}$ and define $h :
\bigcup_{i\in\omega}\{e_i, d_i\} \to k$ such that
$\Phi_A(Z) = kh(e_i) + h(d_i)$ for all $i$. It
follows that no member of $A^*$ extends $h$.

To see that $\{A^* : A\in {\mathcal A}\}$ is a cover of
$\presup{\omega}{k}$ let $f \in \presup{\omega}{k}$. Let $g:\omega \to
k^2$ be defined such that $g(n) = kf(\beta_s(n)) +
f(\beta_g(n))$. Then there is some $A\in \mathcal A$ such that $ g \in
A$. It is easy to check that $f \in A^*$.
\end{proof}
\begin{propo}\label{p:m}
If $ i \geq j$ then $\cov({\mathfrak R}_i) \leq \cov({\mathfrak R}_j)$.
\end{propo}
\begin{proof}
Let $ \bigcup_{\zeta \in \kappa}A_\zeta$ be a
cover of $\presup{\omega}{j}$ by sets in ${\mathfrak R}_j$.
Let $\Phi_\zeta:\PF_j\times  \to \PF_j$ witness that $A_\zeta$ belongs
to  ${\mathfrak R}_i$. 
Define $S:\PF_i \to \PF_j$ by $$S(f)(m) =
\begin{cases}
f(m) &\IF f(m) \in j\\
j-1 &\IF f(m)\notin j\end{cases}$$
and then let $\Psi_\zeta:\PF_i\times  \to \PF_i$ be
defined by $$\Psi_\zeta(f)(m)=
\begin{cases}
\Phi_\zeta(S(f))(m) & \IF m \notin \dom(f)\\
f(m) & \IF m \in \dom(f)\end{cases}$$
Let $B_\zeta = \{f \in \presup{\omega}{i} : (\forall g\in
\PF_i)(\Psi_\zeta(g)\not\subseteq^* f)\}$ and note that if $f \in
\presup{\omega}{i} \setminus \bigcup_{\zeta\in\kappa}B_\zeta$ then $S(f) \in
\presup{\omega}{j} \setminus \bigcup_{\zeta\in\kappa}A_\zeta$.
\end{proof}

\section{Covering Numbers of Many Ros\l onowski Ideals may be Different}

In this section it will be shown 
that any combination of values
for the cardinal invariants $\cov({\mathfrak R}_k)$
is  consistent so long as it does not violate the basic
monotonicity result of  Proposition~\ref{p:m}.

\begin{theor}\label{t:m}
Let $\kappa$ be a  non-increasing 
function from $\omega\setminus 2$ to the uncountable reqular cardinals.
 It is consistent, relative to the consistency of set theory itself, that
$\cov({\mathfrak R}_i) = \kappa(i)$
for each $ i\geq 2$.
\end{theor}

The basic idea of the construction is that a finite support iteration
of length $\kappa(2)$
of countable chain condition partial orders will be constructed. 
At successor stages, Cohen reals will be added and these will be used
to construct trees which will provide an upper bound on
$\cov({\mathfrak R}_i)$. 
At the typical limit stage an approximation to a function witnessing
that $\cov({\mathfrak R}_i)$ is small will have been trapped. A tower
of partial 
functions with respect to $\subseteq^*$ will be constructed and a new
function will be added to the top of this tower. This new function will
prevent the approximation from witnessing
that $\cov({\mathfrak R}_i)$ is small. The countable chain condition
of this tower 
forcing is not an obstacle since this will follow from the genericity
of the construction. More care will have to be taken to preserve the
key property of the trees which guarantee an upper bound on the
covering numbers. The remainder of this section will supply the details.

Let $V$ be a model where there the following hold:
\begin{itemize}
\item $2^{\lambda} \leq \kappa(2)$ for each $\lambda < \kappa(2)$ 
\item There is a $\square_{\kappa(2)}$ sequence --- in other words,
there is family $\{C_\gamma : \gamma\in \kappa(2) \AND \gamma
\text{ is a limit}\}$ such that 
\begin{itemize}
\item each $C_\gamma$ is closed and unbounded in $\gamma$
\item $|C_\gamma| = \cof(\gamma)$ for each
$\gamma$ 
\item if $\delta$ is a limit point of some $C_\gamma$ then
$C_\delta = C_\gamma\cap \delta$
\end{itemize}
\item The following version of $\lozenge$ holds: There is a sequence
$\{D_\alpha\}_{\alpha\in\kappa(2)}$ such that for each $X\subseteq
\kappa(2)$, each closed unbounded set $C\subseteq \kappa(2)$, each
cardinal $\lambda \in \kappa(2)$    and
each $\mu \in \kappa(2)$ there is some $\gamma\in \kappa(2)$ such
that 
\begin{itemize}
\item the order type of $C_\gamma$ is $\lambda$ 
\item  $C_\gamma \subseteq C\setminus \mu $
\item  $D_\zeta = X\cap \zeta$ for each  $\zeta \in  C_\gamma$ which
is a limit of $C_\gamma$. 
\end{itemize}
\end{itemize}
This can be obtained by a strategically closed forcing which is
outlined in the appendix.

The first step is to define a finite support interation of countable
chain condition 
partial orders $\{\Qposet_\alpha\}_{\alpha\in\kappa(2)}$. The
iteration of $\{\Qposet_\alpha\}_{\alpha\in\eta}$ will be denoted by
$\Poset_\eta$. Before proceeding, 
using the cardinal arithmetic hypothesis,
let  all sets of hereditary
cardinality less than $\kappa(2)$ be enumerated by
$\{F_\eta\}_{\eta\in\kappa(2)}$.

If $\alpha = \beta + 2$  then $\Qposet_\alpha$ is simply Cohen
forcing for adding a generic function $c_\alpha:\omega\to\omega$. 
Defined simultaneously with $\Poset_\alpha$ will be trees
$T^{\alpha}_j\subseteq \Omega_j = 
\presup{\kappa(j)}{\kappa(2)}$ and functions $\Theta^{\alpha}_j$
with domain $T^{\alpha}_j$ such that, for each $ j \geq 2$
\begin{itemize}
\item if $\beta \in \alpha$ then  
$T^{\beta}_j \subseteq T^{\alpha}_j$
\item if $\beta \in \alpha$ then  
$\Theta^{\beta}_j \subseteq \Theta^{\alpha}_j$
\item if $\xi \in T_j^\alpha$  then
$1\forces{\Poset_\alpha}{\Theta^\alpha_j(\xi)\in \PF_{j}}$
\item if $\xi$ and $\xi'$ belong to $T_j^\alpha$ and $ \xi\subseteq
\xi'$ then  $1\forces{\Poset_\alpha}{\Theta^\alpha_j(\xi)\subseteq^*
\Theta^\alpha_j(\xi')}$ 
\item if $\xi$ and $\xi'$ are distinct elements of $T_j^\alpha$ of the
same height then  $$1\forces{\Poset_\alpha}{|\{n \in \omega :
\Theta^\alpha_j(\xi)(n)\neq \Theta^\alpha_j(\xi')(n)\}| = \aleph_0}$$
\item if $\alpha$ is a limit then 
 $T^\alpha_j = \bigcup_{\beta\in\alpha} T^\beta_j$
and $\Theta^\alpha_j = \bigcup_{\beta\in\alpha} \Theta^\beta_j$ 
\item if $\alpha = \beta + i$ where $ i \in \{1,2\}$ and $\beta$ is 
 a limit then  $T^\alpha_j = T^\beta_j$
and $\Theta^\alpha_j =  \Theta^\beta_j$ 
\end{itemize}
Notice that by the induction hypothesis, if $F\in \PF_{j}\cup
 \presup{\omega}{{j}}$ and 
 $B_j^\alpha(F)$ is defined to be
 $ \{\xi \in T^{\alpha}_j :
\Theta^{\alpha}_j \subseteq^* F\}$ then $B_j^\alpha(F)$ forms a chain in
$T^{\alpha}_j$. The following additional induction hypothesis will
 play a crucial role in the construction:
\begin{equation}\label{e:kih}
(\forall j \geq 2)(\forall F \in \PF_{j})(|B_j^\alpha(F)| < \kappa(j))
\end{equation}

If $\alpha = \beta + 3$ then let $\varphi(j,\alpha)$ be the least ordinal such
that $F_{\varphi(j,\alpha)} $ is a $\Poset_{\beta + 2}$-name for an element of
$\PF_{{j}}$ which does not appear in the range of
$\Theta^{\beta+2}_j$. (Such an ordinal must exist because $\alpha$ is
a successor and, hence, many new reals have been added at the previous
stage.) Given a generic extension by $\Poset_\alpha$, let 
$F^\alpha_j$ be the interpretation of $F_{\varphi(j,\alpha)}$ in this extension.
Let $\bar{\xi}$ be the lexicographically least member
of $\Omega_j \setminus T^{\alpha}_j$ 
which extends each member of $B_j^{\beta+2}(F^\alpha_j)$ and let $T^\alpha_j =
T^{\beta+1}_j 
\cup \{\bar{\xi}\}$. Note that by \ref{e:kih} the sequence $\bar{\xi}$
belongs to $\Omega_j$. Define 
$\Theta^{\alpha}_j(\bar{\xi})$
by $$\Theta^{\alpha}_j(\bar{\xi})(n) =
\begin{cases}
F^\alpha_j(i) & \IF i \in \dom(F^\alpha_j)\\
c_\alpha(i)  & \IF i \in \omega\setminus \dom(F^\alpha_j) \AND
c_\alpha(i) < j\\
\text{undefined} & \IF
i \in \omega\setminus \dom(F^\alpha_j) \AND
c_\alpha(i) \geq j\end{cases}$$
Notice that this definition will satisfy the induction hypotheses
because of the genericity of $c_\alpha$.
Observe also, that adding a Cohen real does no harm to the induction
hypothesis \ref{e:kih}.

The next step is to
define $\Qposet_\alpha$ when $\alpha$ is a limit or the successor of a
limit ordinal.

\begin{defin}\label{d:tf}
If $\beta$ is an ordinal
and ${\mathcal H} =
\{h_\mu\}_{\mu\in \beta} \subseteq \PF_k$ is such that $h_\mu
\subseteq^* h_\nu$ whenever $\mu\leq \nu$ then the partial order
$\Qposet({\mathcal H})$ is defined to be the set of all functions $f \in
\PF_k$ such that there is some $\mu\in \beta$ such that $f \subseteq^* h_\mu$
ordered under inclusion. If $G$ is a filter on $\Qposet({\mathcal H})$
 then define $f_G = 
\cup G$ and note that if $G$ is a sufficiently generic filter then
$f_G:\omega \to k$. 
\end{defin}
 
Observe  that if  $X\subseteq \beta$ is a cofinal set then
$\Qposet(\{h_\mu\}_{\mu\in X})$ is a dense subset of $\Qposet(\{h_\mu\}_{\mu\in
\beta})A$. This fact will be used in the sequel without further mention.
The function $f_G$ is intended to be used to extend the given chain and
obtain a new partial order extending the given one. However, since
$f_G$ is a total function, it will be necessary
to cut it down  to obtain a member of $\PF_k$. The following
partial order is designed  to do this.

\begin{defin}\label{d:tfa}
If $\Qposet({\mathcal H})$ is as in Definition~\ref{d:tf} and $G$ is a
filter on $\Qposet({\mathcal H})$ then define $$\Sposet(G) = 
\{(a, p) \in [\omega]^{<\aleph_0}\times G : a \cap \dom(p) =
\emptyset\}$$
ordered under coordinatewise inclusion.
If $H \subseteq \Sposet(G)$ is a filter then define 
$A_H = \bigcup_{(a,p)\in H}a$ and define
$f_{G,H} = f_G\restriction (\omega \setminus A_H)$.
\end{defin}

Observe that $\Sposet(G)$ is $\sigma$-centred regardless of the
cofinality of ${\mathcal H}$. Hence $\Qposet({\mathcal H}) *
\Sposet(G)$ has the countable chain condition so long as
$\Qposet({\mathcal H})$ does. Furthemore, $\Qposet({\mathcal H})
\subseteq \Qposet(\{f_{G,H}\})$. The main question to be addressed is:
Do dense sets in $\Qposet({\mathcal H})$ remain dense in
$\Qposet(\{f_{G,H}\})$? The next pair of  lemmas provide some
information on this.

\begin{lemma}\label{l:tfl0}
If ${\mathcal H}\subseteq \PF_k$, $p \in \Qposet({\mathcal H})$, 
$g:l\to k$,
$a \in [\omega \setminus\dom(p)]^{<\aleph_0}$
and $D$ is a dense subset of $\Qposet({\mathcal
H})$
 then there is
$p' \supseteq p$ such that $a\cap \dom(p') = \emptyset$
and
$(p'\restriction(\omega \setminus l)\cup\theta)\cup g
 \in D$
for each $ \theta : a\setminus l \to k$.
\end{lemma}
\begin{proof} This is part of thestandard fusion argument for tree-like forcing.
\end{proof}

\begin{lemma}\label{l:mdps}
Let
${\mathcal H} =
\{h_\mu\}_{\mu\in \beta} \subseteq \PF_k$ be such that $h_\mu
\subseteq^* h_\nu$ whenever $\mu\leq \nu$ and let
 $G$ be $\Qposet({\mathcal H})$-generic over the model $V$. Suppose
also that $H$ is $\Sposet(G)$ generic over $V[G]$. If $D\subseteq
\Qposet({\mathcal H})$ is predense then it remains so in
$\Qposet({\mathcal F})$ for any family ${\mathcal F} \subseteq \PF_k$
such that $f_{G,H} \in \mathcal F$. \end{lemma}
\begin{proof}
>From Lemma~\ref{l:tfl0} it follows that for each dense $D\subseteq
\Qposet({\mathcal H})$ and each $g:l\to k$ the set
 $$D_g = \{(a,p) \in \Sposet(G) :
(\forall \theta (a\setminus l) \to  {k})((p\restriction(\omega \setminus
l)\cup\theta)\cup g  \in D\}$$ is dense in $\Sposet(G)$. Hence, given $f
\supseteq^* f_{G,H}$ choose $l\in\omega$ such that
$f\restriction(\omega \setminus l)
\supseteq f_{G,H}\restriction(\omega \setminus l)$. It may, without
loss of genrality, be assumed that $l\subseteq \dom(f)$  and so it is
possible to let $g = f\restriction l$. Now choose $(a,p ) \in D_g\cap
H$. Let $\theta = f\restriction (a \setminus l)$ and, using the
definition of $D_g$, conclude that $(p\restriction(\omega \setminus
l)\cup\theta)\cup g  \in D$. Since $p\restriction(\omega \setminus
l)\subseteq f_{G,H}\restriction(\omega \setminus l)
\subseteq f$ it follows that $(p\restriction(\omega \setminus
l)\cup\theta)\cup g \subseteq f$ and hence, $f$ extends an element of $ D$.
\end{proof}

Whenever $\alpha $ is a limit ordinal of 
cofinality $\kappa(j)$, the partial
order $\Qposet_\alpha$ will be defined to be of the form
$\Qposet({\mathcal H}_\alpha)$  where
 ${\mathcal H}_\alpha \subseteq \PF_J$ for some $J <\kappa(j)$ is an
increasing tower with respect to $\subseteq^*$ which has cofinality
$\kappa(j)$. 
Moreover, in this case, $\Qposet_{\alpha+1}$ will always be of the form
$\Sposet(G)$ where $G$ is the generic filter
on $\Qposet({\mathcal H}_\alpha)$. Keeping this in mind, let
$H$ be the generic filter on $\Sposet(G)$
and define $H_\alpha = f_{G,H} \in \PF_J$.
The only point which requires elaboration is how to choose ${\mathcal
H}_\alpha$. 

There are three cases to consider. Before proceeding, recall that if
 $C$ is a set of ordinals then $C'$ denotes the Cantor-Bendixon derived
 set of $C$ with respect to the order topology; in other words, $C'$
 is the set of points in $C$ which are limits of $C$. Suppose that for
 each $\xi \in \alpha'$  a family
$\{H^\xi_\gamma\}_{\gamma\in C_\xi} $ has been
 defined.
 To begin, suppose  that the following statement fails:
\begin{equation}\label{e:99}
(\forall \eta \in C_\alpha')
(\forall \bar{\eta} \in C_\eta')
(\forall \xi \in C_{\bar{\eta}})
(H^\eta_\xi = H^{\bar{\eta}}_\xi)
\end{equation}
and there is some $J < j$ such that
 $\{H^\eta_\xi\}_{\xi\in C_\eta} \subseteq \PF_J$
is an increasing tower with respect to $\subseteq^*$
for each $\eta \in C_\alpha$.
In this case let ${\mathcal H}_\alpha$ be any increasing countable
 family; in other words,
 $\Qposet_{\alpha}$ and $\Qposet_{\alpha+1}$ will both be Cohen forcing.
If the statement holds then, for $\xi \in C_\alpha$,
let $H_\xi^\alpha = H_\xi^\eta$ for some (any) $\eta \in
 C_\alpha'\setminus \xi$.
There are two remaining cases.
First, suppose that  $C_\alpha'$ is cofinal in
$\alpha$. In this case 
 ${\mathcal H}_\alpha = \{H^\alpha_\gamma\}_{\gamma
\in C_\alpha}$.
The second case arises if  $C_\alpha'$ is not cofinal in
$\alpha$. Let
$\mu(\alpha)$ be the largest limit of $C_\alpha$ or, if no such limit
exists,  let $\mu(\alpha)= 0$. 
 Suppose also that
$D_\alpha$, as given by the $\lozenge$-sequence, is a 
$\Poset_\alpha$-name and $J < j$
$$1\forces{\Poset_\alpha}{D_\alpha
= \{\Phi^\alpha_\xi\}_{\xi\in\lambda} \AND \Phi^\alpha_\xi:\PF_J \to \PF_J
\text{ witnesses that }\cov({\mathfrak R}_J) \leq \lambda}$$ 
for some $\lambda < \kappa(j)$
Let $\{\gamma_n\}_{n\in\omega}$ enumerate
$C_\alpha\setminus\mu(\alpha)$ in increasing order.
In this case, 
let $H_{\mu(\alpha)}^\alpha = H_{\mu(\alpha)}$ and choose
$H^\alpha_{\gamma_n} $ to be a
$\Poset_\alpha$ name such that  $$1\forces{\Poset_\alpha}{
 H^\alpha_{\gamma_n} =  \Phi^\alpha_{\rho(n)}(H^\alpha_{\gamma_n})}$$
where $\rho(n)$ is the order type of $C_\alpha\cap (\mu(\alpha) + \gamma_n)$. 
Let ${\mathcal H}_\alpha = \{H^\alpha_\eta\}_{\eta \in C_\alpha}$. 

\begin{lemma}\label{l:ccc}
The partial order $\Poset_{\kappa(2)}$ has the countable chain
condition.
\end{lemma}
\begin{proof}
Proceed by induction to show that $$1\forces{\Poset_\alpha}{\Qposet_\alpha\text{ has the countable chain
condition}}$$ for each $\alpha$. 
The countable chain condition for $\Qposet({\mathcal
H})$ is problematic only when  the cofinality of $\beta$ is
uncountable. Indeed, if $\cof(\beta) = \omega$ or $\cof(\beta) = 1$
then $\Qposet({\mathcal H})$ is $\sigma$-centred.  
 If $A\subseteq \Qposet(\{H^\alpha_\gamma\}_{\gamma\in
C_\alpha})$ is a maximal antichain then, using the fact that
$C_\alpha$ is closed and unbounded, it is possible to find some $\gamma \in
C_\alpha$ such that $A\cap \Qposet(\{H^\alpha_\eta\}_{\eta
\in C_\gamma})$
is a maximal antichain. By the induction hypothesis, it follows that
$A\cap \Qposet(\{H^\alpha_\eta\}_{\eta \in C_\gamma})$ is countable. By
Lemma~\ref{l:mdps} it follows that $A\cap \Qposet(\{H^\alpha_\eta\}_{\eta
\in C_\gamma})$ is also maximal in $\Qposet_\alpha =\Qposet(\{H_\eta\}_{\eta
\in C_\alpha}) $.
\end{proof}

Before proceeding some notation will be introduced. 
\begin{defin}
Suppose that
$\Poset \subseteq \Poset'$ and that  $X$ is $\Poset'$ name. The
$\Poset$-name $X\restriction \Poset$ is defined by induction on the
rank of the inductive definition of names. If $X$ is of the form $X
\subseteq \Poset'\times Z$ where $Z$ is a ground model set then
$X\restriction \Poset = X\cap\Poset\times Z$. In general,
$X\restriction \Poset = \{(p,A\restriction \Poset): (p,A) \in X\}$. 
\end{defin}

\begin{lemma}\label{l:lb}
If $G$ is $\Poset_{\kappa(2)}$ generic over $V$ then $\cov({\mathfrak R}_{j})
> \kappa(j)$ in $V[G]$ for $j \geq 2$.
\end{lemma}
\begin{proof}
If $\cov({\mathfrak R}_{j}) \leq \kappa(j)$ then let
$\Phi_\xi:\PF_j \to \PF_j$ be such that $\{\Phi_\xi\}_{\xi\in \lambda}$
 witness this fact for some $\lambda < \kappa(j)$.
Let $\tilde{\Phi}_\xi$ be a name for $\Phi_\xi$ and suppose that $$1\forces{
\Poset_{\kappa(2)}}{\{\tilde{\Phi}_\xi\}_{\xi\in \lambda}\text{ 
 witnesses that }\cov({\mathfrak R}_{j}) \leq \lambda}$$
Let $C$ be a closed unbounded set in $\kappa(2)$ such that for each
$\alpha \in C$  the restricted names
$\tilde{\Phi}_\xi\restriction\Poset_\alpha$ satisfy that
$$1\forces{
\Poset_{\alpha}}{\{\tilde{\Phi}_\xi\restriction \Poset_\alpha\}_{\xi\in \lambda}\text{ 
 witnesses that }\cov({\mathfrak R}_{j}) \leq \lambda}$$
Find some $\gamma$ such that $\cof(\gamma)
=\lambda$, $C_\gamma \subseteq C\setminus \sup(\dom(p))$ and 
$D_\eta = \{\tilde{\Phi}_\xi\restriction
\Poset_\eta\}_{\xi\in\lambda}$ for each $\eta \in 
C_\gamma$. 
 It follows directly from the construction of $\Poset_{\kappa(2)}$ that
$\{H_\rho\}_{\rho \in C_\lambda}$ is an increasing sequence in $
\PF_j$. Moreover, the construction at isolated limit ordinals
guarantees that $H^\gamma_{\rho+1} \supseteq^*
\Phi_\xi(H^\gamma_\rho)$ 
where $\xi$ is the order type of $\rho\cap C_\gamma$ for each
$\rho \in C_\gamma$. This, together with the fact that the
order type of $C_\gamma$ is $\lambda$,  yields that
 $f = f_{G\cap \Qposet(\{H^\gamma_\eta\}_{\eta\in C_\gamma})}$
extends each $\Phi_\xi(H^\gamma_\rho)$ where $\xi$ is the order type
of $\rho\cap C_\gamma$. Hence $f$ does not belong to any of the
members of the ideal ${\mathfrak R}_{j}$ defined by the witnesses $\Phi_\xi$.
\end{proof}
\begin{lemma}\label{l:ub}
If $G$ is $\Poset_{\kappa(2)}$ generic over $V$ then $\cov({\mathfrak R}_{i})
\leq \aleph_{m - j}$ in $V[G]$ provided that $i \geq k_{j}$.
\end{lemma}
\begin{proof}
This follows directly from the induction hypothesis \ref{e:kih}. In
$V[G]$, for each $\alpha \in \kappa(j)$, let $E_\alpha $ be the set
of all $f: \omega \to k_j$ such that there is 
some $\sigma \in T^{\kappa(2)}_{k_j}$ such that the length of $\sigma$
is at least $\alpha$ and $\Theta^{\kappa(2)}_{k_j}(\sigma)\subseteq^* f$.
It is easily verified that $\bigcup_{\alpha\in\kappa(j)} =
\presup{\omega}{k_j}$. The monotonicity esatblished in
Proposition~\ref{p:m} yields the lemma.
\end{proof}

Hence, in order to finish the proof of Theorem~\ref{t:m}, 
it suffices to show that \ref{e:kih} holds.
The first thing to notice is that it suffices to show that the
induction hypothesis holds at a single stage for any particular name
for a function since Cohen genericity will handle the rest. The
 point of the next three lemmas is a stronger version of this assertion
\begin{lemma}\label{l:cgt}
Let $G$ be $\Poset_{\kappa(2)}$ generic over $V$ and $J < j $.
If $\alpha \in \beta \in \kappa(2)$ and $T$ is a $J$-branching
subtree of $\fomom$ 
which belongs to $V[G\cap \Poset_{\alpha}]$ then for any $\xi \in
T_j^\beta\setminus T_j^\alpha$ there are infinitely many integers $i$
such that there is some $i' > i$ so that 
$$\Theta^\beta_j(\xi)\restriction (i' \setminus i) \neq b\restriction
(i' \setminus i)$$ 
for any $b \in \overline{T}$.\end{lemma}
\begin{proof}
Recall that a tree $T$ is said to be $J$-branching of height $n$  if
 $T\subseteq \bigcup_{k\leq n}\presup{k}{\omega}$ and
no node has more than $J$ successors. 
The following fact is easily proved by induction on  $n$: 
If $\{T_i\}_{i\in n}$ is a family of $J$-branching trees 
of height $n$ then $\bigcup_{i\in n}T_i
\not\supseteq \presup{n}{(J+1)}$. A direct corollary of this fact is
that if $T\subseteq \wfomom$ is a $J$-branching tree and $n\in\omega$
then there is a function $f:(i+J^i)\setminus i \to J+1$ such that
$f   \neq b\restriction
((i + J^i) \setminus i)$ 
for any $b \in \overline{T}$. This fact will be used with
Cohen genericity to obtain the desired conclusion.

Before this can be done however, let $T$ and $G$ be given and let $i$
be an arbitrary integer. Let $A$
denote the domain of the interpretation of
$F_{\varphi(j,\beta)}$ in $V[G\cap \Poset_\beta]$.
 Define a tree  $T(i)$  in $V[G\cap
\Poset_\beta]$ by $T(i) =\{t \in T : t\restriction (A\setminus i) \subseteq
F_{\varphi(j,\beta)}\}$ and 
let $\psi_i$ be the order preserving
bijection from
$\omega$ to $\omega\setminus (A\cup i)$. Define $T^*(i) = \{t\circ \psi : t \in
T(i)\}$ and notice that $T^*(i)$ is a $J$-branching tree.
Using this and the Cohen genericity of $c_\beta$
it is possible to apply the observation of the previous paragraph to conclude
 that there are infinitely many
integers $i$ such that $c_\beta\circ\psi\restriction
((i + J^i) \setminus i)   \neq b\restriction
((i + J^i) \setminus i)$ for any $b \in \overline{T^*(i)}$.
Given any such $i$ let $i' = i + J^i +
|A\cap (\psi(J^i)|$. It follows that
$c_\beta\restriction
(i' \setminus i)   \neq b\restriction
(i' \setminus i)$ for any $b \in \overline{T}$. 
\end{proof}

\begin{defin}\label{d:kdl}
If  ${\mathcal H}\subseteq \PF_k$, $g:l\to k$
and $f$ is a $\Qposet({\mathcal H})$-name
such that
$p\forces{\Qposet({\mathcal H})}{f \in \fomom}$ then
a finite subset $a\subseteq \omega$ will be said to  $k$-approximate $f$
with respect to $p$ and $g$  if
\begin{itemize}
\item $\tau : \bigcup_{m \leq |a|}\presup{m\cap a \setminus l}{k} \to \omega$
\item $a \cap \dom(p)=\emptyset$
\item for each $ \theta :a\setminus l \to k$
$$g \cup \theta\cup p\restriction (\omega \setminus l)
\forces{\Qposet({\mathcal H})}{\tau(\theta \restriction j) = f(\theta
\restriction j)}$$
for each $j\leq |a|$.\end{itemize}
Let $f$ be a
$\Qposet({\mathcal H})$-name for a function from $\omega$ to $\omega$
and let $\overline{a}\subseteq a \in
[\omega]^{<\aleph_0}$.
If $G$ is a generic filter on $\Qposet({\mathcal H})$
define $D(f,\overline{a},g)$ to be the set of all
$(a',p) \in \Sposet(G)$ such that $a\subseteq a'$ and
there exists $\overline{a'}\subseteq a'$ such that $\overline{a}$ is a
proper subset of $\overline{a'}$ and
$\overline{a'}$ $k$-approximates $f$ 
with respect to $p$ and $g$.
\end{defin}

It is worth observing that if $a\subseteq \omega$  $k$-approximates $f$
with respect to $p$ and $g$ then the function $\tau$ witnessing this
fact is uniquely defined. Henceforth, this function wil be denoted by
$\tau({a,f,g,p})$.

\begin{lemma}\label{l:tfl}
If  ${\mathcal H}\subseteq \PF_k$, $g:l\to k$
and 
$p\forces{\Qposet({\mathcal H})}{f \in \fomom}$ then,
for any finite subset 
$a\subseteq \omega$ and any $G$ which is a generic filter on
$\Qposet({\mathcal H})$,
$D(f,a,g)$ is dense in
$\Sposet(G)$ below $p$ provided that  
 $a$  $k$-approximates $f$
with respect to $p$ and $g$.
\end{lemma}
\begin{proof} This is a standard  argument based on enumerating
all possible $\theta:a \to k$ and finding a decreasing sequence of
appropriate extensions.\end{proof}

\begin{lemma}\label{l:tf}
If it is given  that 
\begin{itemize}
\item $\cof(\alpha)= \kappa(j)$ 
\item $G$ is $\Poset_{\alpha + 1}$ generic over $V$
\item $f \in \fomom$ in $V[G]$
\end{itemize}
then there is a
$J$-branching tree $T\subseteq \wfomom$ in $V[G\cap \Poset_\alpha]$
 such that ${f \in \overline{T}}$ and $J < j$.\end{lemma} 
\begin{proof} Let
 $\Qposet_\alpha = \Qposet(\{H_\eta\}_{\eta\in\alpha})$. 
Using the countable chain condition of $\Qposet_\alpha$ and the
uncountable cofinality of $\alpha$ it is possible to find a limit
ordinal  $\beta \in C_\alpha$ such that $f$ is a
$\Qposet(\{H_\eta\}_{\eta\in\beta})$-name and 
the name $f$ belongs to $V[G\cap \Poset_\beta]$. Notice that
$\cof(\alpha) = \kappa(j)$ implies that $\{H_\eta\}_{\eta \in
\alpha} \subseteq \PF_J$ for some $J < j$.

Let $H$ be $\Sposet(G\cap \Qposet_\alpha)$ generic over $V[G\cap \Poset_\alpha
*
\Qposet(\{H_\eta\}_{\eta\in\alpha})]$ and let $H_\alpha = f_{G\cap \Qposet_\alpha,H}$.
If $h \in \PF_J$ and $ h \supseteq^* H_\alpha$ 
let $l\in\omega$ be such that  $h\supseteq H_\alpha\restriction(\omega
\setminus l)$ and let $g = h\restriction l$.
Now use Lemma~\ref{l:tfl} to conclude that
there is an infinite chain
$\{a_i^g\}_{i\in\omega}$  such that for each $i$ there is some $p^g_i \in
G$ such that $a_i$ $J$-approximates $f$ with respect to $g$ and $p^g_i$.
Let $A^g= \bigcup_{i\in\omega}a^g_i$ and $\tau^g
=\bigcup_{i\in\omega}\tau({a^g_i,f,g,p^g_i})$. Given $m
\in\omega$  it is possible to extend $h$ to $h'$ such that
$a_m\subseteq \dom(h')$. Let $\theta = h'\restriction a_m$ and
observe that $$g\cup p^g_m\restriction (\omega\setminus a_m) \cup
\theta\forces{\Qposet(\{H_\eta\}_{\eta\in\alpha})}{f(|
\theta\restriction m|)
 = \tau^g_m(\theta\restriction m)}$$
for each $n \leq |a_m|$.
Hence, since $H_\alpha$ extends each $p^g_m$, it follows that
$h'$ forces $f$ to belong to the $j$-branching tree determined by
$\tau^g$.
The desired result now follows directly from
Lemma~\ref{l:mdps}.
\end{proof}

The countable chain condition guarantees that the induction
hypothesis~\ref{e:kih} 
will hold at limit stages of uncountable cofinality, provided that it
holds at all previous stages. The argument at limit stages of
countable cofinality rquires that a bit more care must be taken, but nothing
particular about the forcing is used.

\begin{lemma}\label{l:csih}
The induction hypothesis~\ref{e:kih} holds at limits of countable cofinality,
provided that it holds at all previous stages. 
\end{lemma}
\begin{proof}
Let $\alpha$ have countable cofinality and suppose that
$G$ is $\Poset_\alpha$ generic over $V$. If
$F$ is a function from $\omega$ to $j$ in $V[G]$ then
notice is that, if 
$B^\alpha_j(F)$ has length
$\kappa(j)$ then, by the countable cofinality of
$\alpha$, there is some $\beta \in \alpha$ such that there is a cofinal
subset
 $ B \subseteq B^\beta_j(F)$. This determines
the  branch through $T_j^\beta$ in $V[G\cap \Poset_\beta]$.
 Hence, it
suffices to show that if $B\subseteq T_j^\beta$ is a branch of length
$\kappa(j)$ in $V[G\cap \Poset_\beta]$ and $F$ is in $V[G\cap
\Poset_\alpha]$ then $B \not\subseteq B^\beta_j(F)$.

To this end, let $B$ be a $\Poset_\beta$ name for a long branch
through $T_j^\beta$ and $F$ a $\Poset_\alpha$-name. 
Let $\{\beta_n\}_{n\in\omega}$ be a sequence of ordinals cofinal in
$\alpha$ such that $\beta_n > \beta$ for each $n$.
For any $p \in \Poset_\alpha$ define $F_p = \{(i,j) :
p\forces{\Poset_\alpha}{F(i) = j}\}$.  
It will first be shown that for each $n\in\omega$ the set
$$D(n) = \{ q \in \Poset_{\beta_n} : (\exists \sigma \in B)(\forall
r \leq q)( F_r \not\supseteq^* \Theta^\beta_j(\sigma))\}$$
is dense in $\Poset_{\beta_n}$. To see that this is so, suppose that $
q\in \Poset_{\beta_n}$ is such that for each $\sigma \in B$ and
$\bar{q} \leq q$ there is some $ r \leq \bar{q}$ such that $F_r
\supseteq^* \Theta^\beta_j(\sigma)$. 
Then let $\bar{F}$ to be the $\Poset_{\beta_n}$-name defined by
$p\forces{\Poset_{\beta_n}}{\bar{F}(i)=j}$ if and only if
$p\forces{\Poset_{\alpha}}{{F}(i)=j}$. It follows that
$q\forces{\Poset_{\beta_n}}{\bar{F}\supseteq^*\Theta^\beta_j(\sigma)}$
for each $\sigma\in B$ contradicting the induction hypothesis.

Using the density of each $D(n)$, let ${\mathcal A}_n\subseteq D(n)$
be a maximal antichain and, for each $q \in {\mathcal A}_n$, let
$\sigma_q^n \in B$ witness that $q\in D(n)$. Let $\sigma \in B$ be such
that $\sigma \supseteq \sigma^n_q$ for each ${n\in\omega}$ and $q \in
{\mathcal A}_n$. Now suppose that $p \in
\Poset_\alpha$ is such that $p\forces{\Poset_\alpha}{F\cup
(\Theta^\beta_j(\sigma)\restriction m)
\supseteq\Theta^\beta_j(\sigma)}$.  Let $n$ be such that $p \in
\Poset_{\beta_n}$ and choose $q \in {\mathcal A}_n$ such that there is
some $ r \in \Poset_{\beta_n}$ such that $r \leq q$ and $r \leq p$. 
Since $q\in D(n)$ it follows that $F_r\not\supseteq^*
\Theta^\beta_j(\sigma^n_q) \subseteq^* \Theta^\beta_j(\sigma)$. Hence,
there is some $i > m$ in the domain of $\Theta^\beta_j(\sigma)$ such
that either 
$r\forces{\Poset_\alpha}{F(i)\neq
\Theta^\beta_j(\sigma)(i)}$ or $r$ does not decide a value for
$F(i)$. The first case directly contradicts that $r \leq p$ and, 
in the second case, it is possible to  extend $r$ to $r'$
such that $r'\forces{\Poset_\alpha}{F(i)\neq
\Theta^\beta_j(\sigma)(i)}$. This agina yields a contradiction.
\end{proof}

It remains to consider successor ordinals. If $\alpha = \beta + 1$ and
$\beta$ itself is a successor, then $\Qposet_\alpha$ is
$\sigma$-centred and, hence, a standard argument shows that it
preserves the induction hypothesis. If $\beta$ is a limit of countable
confinality, then $\Qposet_\alpha$ is also $\sigma$-centred.
 So the only problem may arise whent $\beta$ is a limit of
uncountable cofinality. 

\begin{lemma}\label{l:ucih}
Suppose that
  $\alpha$ is a limit ordinal of uncountable cofinality.
 Given  that each preceding stage
 satisfies the induction hypothesis~\ref{e:kih}, the partial order
  $\Poset_{\alpha 
 + 1}$ will also satisfy the induction hypothesis.
\end{lemma}
\begin{proof}
Let $G$ be $\Poset_\alpha$ generic over $V$ and argue in $V[G]$.
There are two types of branches which might provide difficulties.
 To begin, consider
 branches which occur at some stage before $\alpha$.
Let $B$ be a branch through $T^\beta_j$ of length $\kappa(j)$ in $V[G\cap
\Poset_\beta]$ and let $F$ be a $\Qposet_\alpha$-name for a function 
from $\omega$ to $j$ such that
$$1\forces{\Qposet_\alpha}{(\forall \sigma \in
B)(F\supseteq^* \Theta^\beta_j(\sigma))}$$ 
 If $\kappa(j) > \cof(\alpha)$ then $\Qposet_\alpha$ has a dense
subset of cardinality $\cof(\alpha)$ and a pigeonhole argument shows
that there is some $M\in\omega$ and a single condition $q \in
\Qposet_\alpha$ such that the set of $\sigma \in B$ such that
$q\forces{\Qposet_\alpha}{F\cup\Theta^\beta_j(\sigma)\restriction M\supseteq
\Theta^\beta_j(\sigma)}$ is cofinal in $B$. On the other hand,
 if $\kappa(j) < \cof(\alpha)$ then $\Qposet_\alpha$ has
$\kappa(j)$ as a precalibre. In this case it is possible to find
$\{q_\sigma\}_{\sigma \in B'}$ a centred subset of $ \Qposet_\alpha$
and $M \in \omega$  
such that
$q_\sigma\forces{\Qposet_\alpha}{F\cup\Theta^\beta_j(\sigma)\restriction M\supseteq
\Theta^\beta_j(\sigma)}$
 for each $\sigma \in B'$ and, furthermore, $B'$ is a cofinal subset of
$B$. In either case a contradiction is obtained since it follows that
$\Poset_\alpha$ violates the induction hypothesis.
 Hence, it may be assumed
that $\kappa(j) = \cof(\alpha)$. 
Using the countable chain condition of $\Qposet_\alpha$, let $\xi\in
C_\alpha\setminus \beta$ be 
such that  $F\restriction
\Qposet(\{h_\eta\}_{\eta\in \xi})$ is a $\Qposet(\{h_\eta\}_{\eta\in
\xi})$-name.  
Since the cardinality of $\Qposet(\{h_\eta\}_{\eta\in \xi})$ is less
than that of $B$, it follows that $1\forces{\Qposet(\{h_\eta\}_{\eta\in
\xi})}{F\not\supseteq^* \Theta^\beta_j(\sigma))}$ for some fixed
$\sigma \in B$. Now use Lemma~\ref{l:mdps} to conclude that the dense
sets witnessing this remain dense in $\Qposet_\alpha$.

The second possibility is that  a cofinal branch is added to
$T^\alpha_j$. To see that this can not happen,
suppose that $1\forces{\Qposet_\alpha}{F:\omega \to j}$
Then, by Lemma~\ref{l:tf}, there is some $J$-branching tree $T$ such
that $J < j$ and  $1\forces{\Qposet_\alpha}{F \in \overline{T}}$. Since
$\alpha$ has uncountable cofinality and the iterands all have
the countable chain condition, it follows that 
if $G$ is
a generic set for $\Poset_{\kappa(2)}$ then there is some $\beta 
\in \alpha$ such that $T$ belongs to $V[G\cap \Poset_\beta]$. 
Choose $\sigma\in B^\alpha_j(F)\setminus T^\beta_j$.
Now use
Lemma~\ref{l:cgt} to 
obtain a contradiction.
\end{proof}
\section{Appendix}
A brief note regarding the consistency of the required combination of
$\square$ and $\lozenge$ may be helpful to some readers.
To obtain the required initial model, begin with a model where $\lambda$ is
regular and $2^\kappa \leq \lambda$ for $\kappa < \lambda$. Let
$\Poset$ be the partial order consisting of initial segments of the
required $\square$ and $\lozenge$ sequence. To be precise, $p \in
\Poset$ if and only if $p$ is a function defined on some
$\alpha \in \lambda$ such that
\begin{itemize}
\item $p(\eta) = (C_\eta,D_\eta)$ for each $ \eta \in \alpha$
\item $D_\eta \subseteq \eta$
\item $C_\eta \subseteq \eta$ is closed and unbounded in $\eta$
\item if $\eta \in \alpha$ and $ \xi\in C_\eta'$ then $C_\xi =
C_\eta\cap \xi$.
\end{itemize}
This partial order has size $\lambda$ and is stategically
$\lambda$-closed.

\end{document}